\documentclass[11pt]{article}
\usepackage{amsfonts}
\usepackage{amssymb}
\newtheorem{thm}{Theorem}
\newtheorem{claim}{Claim}
\newtheorem{prop}{Proposition}
\newtheorem{dfn}{Definition}
\newtheorem{cor}{Corollary}

\newtheorem{Lemma}{Lemma}
\newcommand{\bpr}{{\bf Proof. }}
\newcommand{\qed}{$\Box$}
\def\g{\mathfrak g}
\def\sl{\mathfrak {sl}}
\def\spin{\mathfrak {spin}}
\def\a{\mathfrak {a}}
\def\sp{\mathfrak{sp}}
\def\so{\mathfrak{so}}

\def\h{\mathfrak h}

\def\z{\mathfrak z}
\def\s{\mathfrak s}
\def\dim{\hbox{\rm dim}\,}
\def\codim{\hbox{\rm codim}\,}

\title{Orthogonal linear group-subgroup pairs with the same invariants}
\author{S.Solomon\thanks{Hebrew University of Jerusalem, Israel. E-mail: gsasha@math.huji.ac.il}}
\input xy
\xyoption{all}
\begin{document}
\maketitle
\begin{abstract}

The main theorem of Galois theory states that there are no finite
group-subgroup pairs with the same invariants. On the other hand, if we
consider complex linear reductive groups instead of finite groups, the
analogous statement is no longer true: There exist counterexample
group-subgroup pairs with the same invariants.  However, it's possible to
classify all these counterexamples for certain types of groups. In
\cite{irrex}, we provided the classification for connected complex
irreducible groups, and, in this paper, for connected complex orthogonal
groups, i.e., groups that preserve some non-degenerate quadratic form.
                                                                                
\end{abstract}

\section{Introduction}

\subsection{Setting the problem}
Let $G$ be an algebraic group acting on an irreducible algebraic
variety $X.$ Consider the action of $G$ on the
field of the rational functions $k(X).$  Let $k(X)^G$ denote the subfield of rational $G$-invariants.
Then $G$ is a group of automorphisms of the field extension $k(X)$
over $k(X)^G,$ that is, $G\subseteq Aut(k(X),k(X)^G).$

Suppose $G$ is a finite group. Then $k(X)/k(X)^G$ is a Galois extension. The main theorem of
the Galois theory states that $G=Aut(k(X),k(X)^G).$  Equivalently, any proper subgroup $H\subsetneq G$
corresponds to a nontrivial extension $k(X)^H\supsetneq k(X)^G.$ Simply saying, $H$ always has more invariants than
$G.$

In other words, a finite group action is uniquely determined by its invariants.  A natural question
arises if this is true for other classes of groups.
The answer in general is "no". The simplest counterexample is
$Sp(V)\subset SL(V),$ where V is an even dimensional vector space.
For both the group and the subgroup,  the only invariants  are constants.

However, it is possible to describe all the counterexamples for certain types
of group actions.
\begin{dfn} We call a triple $(H,G,X)$ exceptional, if the fields of rational invariants of $H$ and $G$ coincide:  $k(X)^H =~k(X)^G.$ If
$H=G$ we say that the triple $(H,G,X)$ is trivial.  \end{dfn}

A general problem is to classify exceptional triples. We just saw that any exceptional triple $(H,G,X)$ with $G$
finite is trivial. In
\cite{irrex}, we classified exceptional triples $(H,G,V),$ where $V$ is a complex (finite dimensional) vector space, $H,G$ are connected irreducible  linear groups, see section~\ref{irrsummary} here for a partial summary. In this paper, we allow the groups
to be reducible but require that they preserve a non-degenerate quadratic form.

\subsection{Background}
To the best of our knowledge, this problem has not been considered before even for reductive linear groups.
However, there are two important classical results in this direction.

E.B.Vinberg in 1959 started, Sato and Kimura \cite{SK} in 1977, and Shpiz\cite{Shp} in 1978, completed the
classification of irreducible complex linear groups $G\subset GL(V)$ acting with an orbit open in $V,$
so-called locally transitive groups (alternatively, $V$ is called a  prehomogeneous vector space). Such
groups cannot have non-trivial polynomial or rational invariants. Thus, any locally transitive group $G$ with a
locally transitive subgroup $H$ make an  exceptional triple $(H,G,V).$

D.Montgomery and H.Samelson\cite{MS} in 1943, and A.Borel\cite{Bor} in 1950, described real  groups that act
transitively on spheres.  Any triple $(H,G,S)$ of a group $G$ and a subgroup $H$ transitive on a sphere $S$ is exceptional.

\subsection{Geometric meaning}
By Rosenlicht's theorem, a triple
$(H,G, X)$  is exceptional if and only if $H$ is locally transitive on a $G$-orbit in general position $O_X(G),$ that is, the closures of orbits in general position of $H$ and
$G$ coincide: $\overline{O_X(H)}=\overline{ O_X(G)}.$ Note that for stable actions, i.e., where the orbits in general position are closed,  it means $O_X(H)=O_X(G).$ In particular, this holds for orthogonal groups $H$ and $G.$

As a consequence, we obtain

\begin{Lemma}\label{ratreg} Suppose $X$ is irreducible quasiaffine, and the polynomial $G$- and $H$-invariants  separate (respective) orbits in general position.  Then $(H,G, X)$ is an exceptional triple if and only if the polynomial invariants of $H$ and $G$ coincide: $k[X]^G=k[X]^H.$ In particular, this holds in case $H,G$ are semisimple, or in case $X=V$ is a vector space, and $H,G$ are orthogonal.
\end{Lemma}
\bpr "If" is trivial. Suppose $(H,G, X)$ is exceptional. By \cite{55} Proposition 3.4, we have
$tr.deg$ $k[X]^G=\codim O_X(G)=\codim O_X(H)=tr.deg$ $k[X]^H,$ hence $k[X]^G=k[X]^H.$ "In particular" statement follows from \cite{55} Theorem 3.3\qed

 \subsection{Main results}
 Let $X=V$ be a vector space over $\mathbb C.$
The main result is the following Theorem~\ref{ortthm} that classifies exceptional triples $(H,G,V)$ under the assumption that $H,G$ are connected semisimple orthogonal groups. Furthermore, Proposition~\ref{reduc} gives a criterion for an exceptional triple $(H,G,V),$ where  $H,G$ are connected reductive orthogonal and at least one of them is not semisimple.

Assume $H\subset G\subset GL(V)$ are semisimple linear groups.

We assign a three level diagram to the triple $(H,G,V).$
The upper and the middle level
vertices correspond to normal subgroups of $H$ and $G,$ respectively.  If a normal subgroup of $H$ projects
nontrivially on a normal subgroup of $G,$ we connect the vertices with an edge, and sometimes specify the projection on the side of the edge. A double edge means that the projection
is an isomorphism. We substitute one of the ends of a double edge with a circle.
The lower level vertices
correspond to a decomposition of $V$ into $G$-submodules. We connect each submodule vertex with all normal
subgroups of $G$ that act nontrivially on this submodule. Unspecified subgroups and submodules, that is,
circles, can be considered arbitrary.

A diagram that has two connected components corresponding to triples $(H_1,G_1,V_1)$ and $(H_2, G_2, V_2)$ encodes the direct sum of these triples, i.e., the triple  $(H_1\times H_2,
G_1\times G_2, V_1\oplus V_2).$ The direct sum of exceptional triples is exceptional, and,  vice versa, if the direct sum is exceptional, then every summand is exceptional. A triple (or a linear group) that cannot be decomposed into a nontrivial direct sum  is called {\it indecomposable}.

We say a triple $(H,G,V)$ is {\it locally trivial} if the restriction $(H,G,U)$ to every proper $G$-submodule $U\subsetneq V$ is a trivial triple.
Also, we call
$G$  {\it strongly faithful} if all nonzero $G$-submodules of $V$ are faithful. We
say a triple $(H,G,V)$ is strongly faithful if $G$ is such.

Finally, we may always assume for our purposes that $H$ is maximal in $G.$ Indeed, if $(H,F,X)$ and $(F,G,X)$ are exceptional triples then $(H,G,X)$ is also exceptional. Conversely, if $(H,G,X)$ is exceptional, then $(F,G,X)$
is exceptional for any $F\supset H.$

\begin{thm}\label{ortthm}
 Let $(H,G,V)$ be an indecomposable nontrivial exceptional triple, where $H,$ $G$ are connected semisimple orthogonal, and $H$ is maximal in $G.$
Then $(H,G,V)$ is either  locally trivial with diagram one of $\bf{T1-T9},$  or strongly faithful with
diagram one of $\bf{F1-F10},$  or $(H,G,V)$ has a diagram $\bf{O},$ where

 (i) in diagram $\bf{T9}$ one of the two indecomposable components of group $G$ satisfies the explicit criterion of Lemma~\ref{yak}; and

(ii) in diagram $\bf{O}$ the action of the subgroup $SO_8\subset G$ on $V_1\oplus V_2$ is isomorphic, up to a trivial component, to $t\phi_1,$  $t\le 3,$ where $\phi_1$ stands for the tautological action of $SO_8$ on $\mathbb C^8.$

Conversely, all the abovementioned diagrams correspond to exceptional triples.

\end{thm}

In the diagrams below, $\nu$ stands for $\nu(H,G,V):=$ tr.deg.$\mathbb C
(V)^H =$ tr.deg.$\mathbb C(V)^G,$ $\nu=4+$ tr.deg.$\mathbb C(U)^{\check G}$ in diagram {\bf T8}. See section~\ref{notation} for other notation.
\par\bigskip
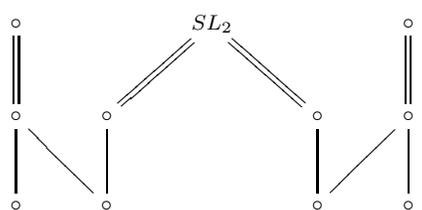
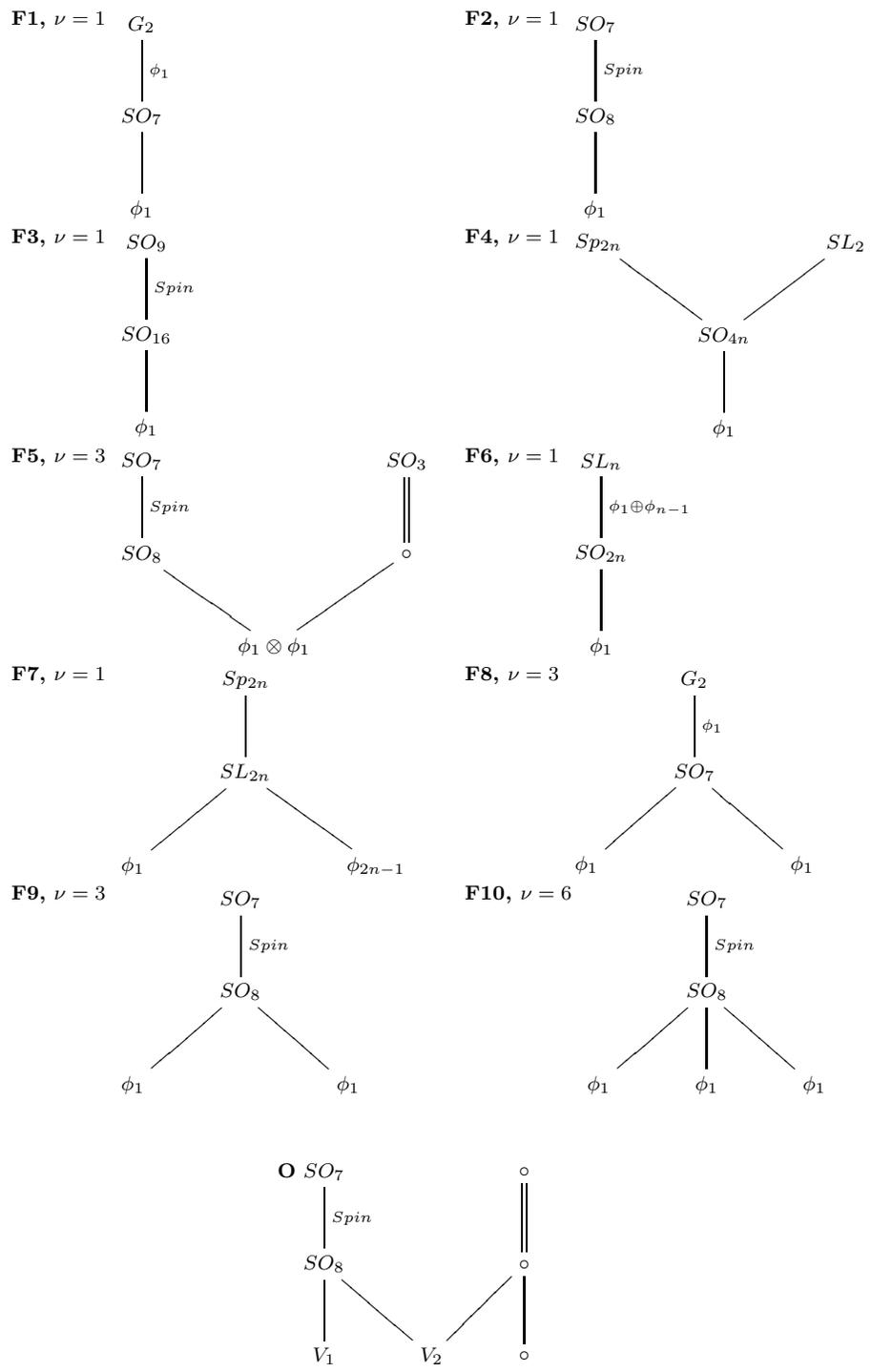
\begin{scriptsize}

\begin{tabular}{ll}

\bf{T1}, $\nu=2$ \xymatrix{
& &{SL_4} \ar@{=}[dl]\ar@{=}[dr]& \\
 &{\circ} \ar@{-}[dl]\ar@{-}[dr] & &{\circ} \ar@{-}[d]\\
{\phi_1}& &{\phi_3}& {\phi_2} \\
}

&

\bf{T2}, $\nu=2$ \xymatrix{
 &{SO_7} \ar@{=}[dl]\ar@{=}[dr]& \\
 {\circ} \ar@{-}[d]& & {\circ} \ar@{-}[d]\\
{\phi_1}& &{\phi_3}\\
}
\\

\bf{T3}, $\nu=2$ \xymatrix{
 &{SO_8} \ar@{=}[dl]\ar@{=}[dr]& \\
 {\circ} \ar@{-}[d]& & {\circ} \ar@{-}[d]\\
{\phi_1}& &{\phi_3}\\
}

&

\bf{T4}, $\nu=2$ \xymatrix{
 &{SO_{12}} \ar@{=}[dl]\ar@{=}[dr]& \\
 {\circ} \ar@{-}[d]& & {\circ} \ar@{-}[d]\\
{\phi_1}& &{\phi_5}\\
}
\\

\bf{T5}, $\nu=4$ \xymatrix{
& &{SO_7} \ar@{=}[dl]\ar@{=}[dr]& \\
 &{\circ} \ar@{-}[dl]\ar@{-}[dr] & &{\circ} \ar@{-}[d]\\
{\phi_1}& &{\phi_1}& {\phi_3} \\
}

&

\bf{T6}, $\nu=4$ \xymatrix{
& &{SO_8} \ar@{=}[dl]\ar@{=}[dr]& \\
 &{\circ} \ar@{-}[dl]\ar@{-}[dr] & &{\circ} \ar@{-}[d]\\
{\phi_1}& &{\phi_1}& {\phi_3} \\
}
\\

\multicolumn{2}{l}{
\bf{T7}, $\nu=7$ \xymatrix{
& &{SO_8} \ar@{=}[dl]\ar@{=}[dr]& \\
&{\circ} \ar@{-}[dl]\ar@{-}[d]\ar@{-}[dr]& &{\circ} \ar@{-}[d]\\
{\phi_1}&{\phi_1}&{\phi_1}&{\phi_3} \\
}
}  \\

\multicolumn{2}{l}{
{\bf T8} \xymatrix{
 &{SO_8} \ar@{=}[dl]\ar@{=}[dr]& & & {SO_3} \ar@{=}[d]& {\check G} \ar@{=}[d]\\
{\circ} \ar@{-}[d]& &{\circ} \ar@{-}[dr] & & {\circ} \ar@{-}[dl]\ar@{-}[dr]& {\circ} \ar@{-}[d]\\
{\phi_3}& & &{\phi_1\otimes \phi_1}& &{U}\\
}
}  \\

\multicolumn{2}{l}{
\bf{T9}
\xymatrix{
 {\circ} \ar@{=}[d]&  &{SL_2}\ar@{=}[dl]\ar@{=}[dr]& &{\circ} \ar@{=}[d]\\
{\circ} \ar@{-}[dr]\ar@{-}[d] &{\circ} \ar@{-}[d]& &{\circ}\ar@{-}[d]&{\circ}\ar@{-}[dl]\ar@{-}[d]\\
{\circ}&{\circ}&  & {\circ}& {\circ}\\
}
}  \\
\end{tabular}
\par\smallskip

\begin{tabular}{ll}
\bf{F1}, $\nu=1$ \xymatrix{
 {G_2} \ar@{-}[d]^{\phi_1}\\
 {SO_7}\ar@{-}[d]\\
 {\phi_1}\\
}
&

\bf{F2}, $\nu=1$ \xymatrix{
 {SO_7} \ar@{-}[d]^{Spin}\\
 {SO_8}\ar@{-}[d]\\
 {\phi_1}\\
}
\\

\bf{F3}, $\nu=1$ \xymatrix{
 {SO_9} \ar@{-}[d]^{Spin}\\
 {SO_{16}}\ar@{-}[d]\\
 {\phi_1}\\
}
&

\bf{F4},  $\nu=1$ \xymatrix{
 {Sp_{2n}} \ar@{-}[dr]& &{SL_2} \ar@{-}[dl]\\
 &{SO_{4n}} \ar@{-}[d]& \\
 &{\phi_1}& \\
}
\\

\bf{F5}, $\nu=3$ \xymatrix{
 {SO_7}\ar@{-}[d]^{Spin}& &{SO_3} \ar@{=}[d]\\
 {SO_8} \ar@{-}[dr]& &{\circ} \ar@{-}[dl]\\
 &{\phi_1\otimes \phi_1}& \\
}
&

\bf{F6}, $\nu=1$ \xymatrix{
 {SL_n} \ar@{-}[d]^{\phi_1\oplus\phi_{n-1}}\\
 {SO_{2n}}\ar@{-}[d]\\
 {\phi_1}\\
}
\\

\bf{F7}, $\nu=1$ \xymatrix{
 &{Sp_{2n}} \ar@{-}[d]&\\
 &{SL_{2n}}\ar@{-}[dl]\ar@{-}[dr]&\\
 {\phi_1}& &{\phi_{2n-1}}\\
}
&

\bf{F8}, $\nu=3$ \xymatrix{
 &{G_2} \ar@{-}[d]^{\phi_1}&\\
 &{SO_7}\ar@{-}[dl]\ar@{-}[dr]&\\
 {\phi_1}& &{\phi_1}\\
}
\\

\bf{F9}, $\nu=3$ \xymatrix{
 &{SO_7} \ar@{-}[d]^{Spin}&\\
 &{SO_8}\ar@{-}[dl]\ar@{-}[dr]&\\
 {\phi_1}& &{\phi_1}\\
}
&

\bf{F10}, $\nu=6$ \xymatrix{
 &{SO_7} \ar@{-}[d]^{Spin}&\\
 &{SO_8}\ar@{-}[dl]\ar@{-}[d]\ar@{-}[dr]&\\
 {\phi_1}& {\phi_1}&{\phi_1}\\
}
\\

\\
\end{tabular}

\[\bf{O}\xymatrix{
{SO_7} \ar@{-}[d]^{Spin}& &{\circ}\ar@{=}[d]\\
{SO_8} \ar@{-}[d]\ar@{-}[dr]& &{\circ} \ar@{-}[d]\ar@{-}[dl]\\
{V_1}&{V_2}&{\circ}\\
}\]

\end{scriptsize}

\begin{prop}\label{reduc}
Let $H\subsetneq G\subset SO(V)$ be connected reductive orthogonal groups. Let $Z(G)$ and $G'$ denote the center and the commutator subgroup of $G,$ and let $T\subset Z(G)$ denote a complement subgroup to the projection of $Z(H)$ on ${Z(G)}.$  Then $(H,G,V)$ is exceptional if and only if

(i) the triple $(H',G',V)$ is exceptional, and

(ii) $T$ acts trivially on $\mathbb C[V]^{H'}=\mathbb C[V]^{G'}.$

\end{prop}

\section{Acknowledgments} No words would suffice to express my gratitude to E.B.Vinberg, my teacher.  I am grateful to A.Shalev, my
advisor at Hebrew University of Jerusalem, Israel. I would like also to thank Michael Larsen for a helpful
discussion, D.I.Panyushev for remarks, and O.Yakimova for pointing out Lemma~\ref{yak}. My family provided
me with continuous inspiration.

\section{Notation and agreements}\label{notation}

All groups $G,H,F,\dots$ are connected complex reductive algebraic groups, $X,Y,\dots$ stand for complex algebraic varieties, and  $V,U,W,\dots$ denote finite dimensional complex vector spaces, unless specified otherwise.

For $G$ acting on an irreducible $X,$
$\mathbb C(X)^G$ and $\mathbb C[X]^G$  denote  the field of {\it rational invariants}, and the algebra of
{\it polynomial invariants}.

For an ex\-cep\-tional triple $(H,G,X),$ $\nu (H,G,X):=$ tr.deg.$\mathbb C
(X)^H =$ tr.deg.$\mathbb C(X)^G.$

Let $G$ be a simply connected semisimple group, and $H\subset G$ be a maximal semisimple subgroup. Then
either (1) $G=G_1\times G_2\supset H=H_1\times G_2$, where $G_1$ is simple, and $H_1\subset G_1$ is a maximal
semisimple subgroup, or (2) $H=H_{12}\times F\subset G=H_1\times H_2\times F,$ where $H_{12}\cong
H_1\cong H_2$ are simple, and $H_{12}$ is embedded diagonally in $H_1\times H_2$ (\cite{Dyn}).

\begin{dfn}\label{maximal}
We say $H$ is a {\it straight} subgroup if (1) holds, and a {\it diagonal} subgroup if (2) holds. We use the
same terminology for the tangent algebras $\h\subset\g.$
\end{dfn}

$O_X(G)$ is an {\it orbit in general position} of $G,$ $G_*(X)$ is a {\it stabilizer in general position (SGP)} of $G.$

$Z(G)$ is the {\it center}, and $G'$ is the {\it commutator subgroup} of $G.$

$X/G:=$Spec$\mathbb C [X]^G$ denotes {\it the categorical quotient}.

A triple $(H,G, X)$ is called {\it semisimple, etc.,} if  both $H$ and $G$ are semisimple, etc.
A triple $(H,G, V)$ is called {\it irreducible, orthogonal, etc.}, if both $H$ and $G$ are irreducible, orthogonal, etc., on $V.$

For a simple, connected, simply connected group $G$ of rank $r$ denote by $\phi_1,\dots,\phi_r$ the fundamental
representations of $G$ in standard ordering (see, for example \cite{Ela1}). Sometimes we denote by $Spin$
the spin representation $\phi_k$ or half-spin representations $\phi_k$ and $\phi_{k-1}$ of $B_k$ or $D_k$.

\section{Where is the difficulty}\label{diff}

\begin{Lemma}\label{factor} Let $(H,G, X)$ be an exceptional triple. Suppose $G$ acts on a variety
$Y,$ and $\pi:X\longrightarrow Y$ is a surjective $G$-equivariant homomorphism.  Then $(H,G,Y)$
is also an exceptional pair.
\end{Lemma}

\bpr Consider $\pi^*:  \mathbb C(Y)\longrightarrow \mathbb C(X),$ where $(\pi^*f)(v)=f(\pi (v)).$ Then
$\pi^*$ is $G$-equivariant and injective.  Hence, $\mathbb C(Y)^H\ne\mathbb C(Y)^G$ would imply $\mathbb
C(X)^H\ne\mathbb C(X)^G.$\qed

Assume $X=V$ is a vector space, and $H,G$ are reductive.

\begin{cor}\label{sub}   Let $U\subset V$ be a $G$-submodule.  If $(H,G,V)$ is exceptional then $(H,G,U)$  is also  exceptional.  \end{cor}

\bpr Since $G$ is reductive, there exists a $G$-invariant subspace $U'$ such that $V=U\oplus U'.$ Denote by
$\pi:  V\longrightarrow U$ the projection on $U$ parallel to $U'.$ Then $\pi$ is $G$-equivariant, and, by
Lemma~\ref{factor}, $(H,G,U)$ is exceptional.\qed

Since $V$ decomposes into a direct sum of $G$-submodules, the triple $(H,G, V)$ is in a certain sense a sum of
exceptional triples, in which the group is irreducible. Note, however, that if
$(H,G,U)$ and $(H,G,W)$ are exceptional triples, the triple $(H,G,U\oplus W)$ is not at all guaranteed to be exceptional. For example,
the triple $(Sp_{2n}, SL_{2n},\mathbb C^{2n})$ is exceptional, and  the triple $(Sp_{2n}, SL_{2n},\mathbb C^{2n}\oplus\mathbb C^{2n})$ is not.

Our main two questions therefore are:

1. What are the exceptional triples $(H,G, V)$ with $G$ irreducible, and

2. How exactly two exceptional triples can combine to make another exceptional triple.

In paper \cite{irrex}, we partially answered the first question, namely, classified  exceptional triples $(H,G,V)$
with both $H$ and $G$ irreducible on $V$ (for convenience, we summarized it here in section~\ref{irrsummary}). The
second question presents a much harder problem due to a generally complicated structure of $V$ as a
$G$-module. The main reason is that  proper submodules of $V$ do not have to be faithful even if
$V$ is faithful. This phenomenon is sometimes called ``blinking
kernels''. As a simplest example, consider a group $G:=G_{12}\times G_{23}\times G_{13}$ acting on a vector space $V=
V_1\oplus V_2\oplus V_3,$ where $G_{ij}$ acts nontrivially on
$V_i$ and $V_j,$ and trivially on the other summand.
In case of blinking kernels, it's very hard to see the orbit structure of $V.$

\subsection{Branching law}

To approach our classification, we need to know how an irreducible $G$-module branches  into irreducible
$H$-modules.  Assume $(H,G, V)$ is an exceptional triple, $G$ is orthogonal. It turns out that for such triples the
branching law is simple: an irreducible $G$-module decomposes into a sum of at most
two $H$-irreducible modules. We describe this law precisely in Lemmas~\ref{submodules} and~\ref{spacedecom}.

Denote by $(\cdot,\cdot)$ an invariant symmetric bilinear form on $V$. A subspace $W\subset V$ is called
non-degenerate, if the restriction of $(\cdot,\cdot)$ on $W$ is non-degenerate.

\begin{Lemma}\label{submodules}
Suppose $W\subset V$ is a non-degenerate $H$-submodule. Then $W$ is also a $G$-submodule. Consequently, if
$H$-action on $W$ is orthogonal, then $W$ is a $G$-submodule.
\end{Lemma}

\bpr Denote by $W^{\perp}$ the orthogonal complement of $W$ in $V.$ Since $W$ is non-degenerate, we have
$V=W\oplus W^{\perp}.$ Take the orthogonal projection $P: V\longrightarrow W^{\perp}.$ Since $W$ is
$H$-invariant, $P$ commutes with $H$. Denote $F(v):=(Pv, Pv).$ $F$ is an $H$-invariant, and,
therefore, a $G$-invariant polynomial. Since $W=Ker$ $F,$ $W$ is $G$-invariant. \qed

As known from linear algebra, a minimal orthogonal $H$-submodule $W\subset V$ is either irreducible, or
$W=U\oplus U^*$, where $U$ is irreducible and non-orthogonal.

\begin{dfn}\label{I-III} We say that a minimal orthogonal $G$-submodule $W\subset V$ is solid if  $W$ is $H$- and $G$-irreducible;
half-split if  $W$ is $G$-irreducible and $W=U\oplus U^*,$  $U$ is $H$-irreducible; and split if $W=U\oplus
U^*,$  $U$ is $H$- and $G$-irreducible. We will also say that a submodule $W$ is solid (half-split, split), if all of its minimal orthogonal $G$-submodules are such.
\end{dfn}

Lemma~\ref{submodules} implies

\begin{Lemma}\label{spacedecom}
$V$ decomposes into direct sum of solid, half-split and split $G$-submodules.
\end{Lemma}






\section{Preliminaries} This is a collection of facts we use throughout the paper.

Let $G_x\subset G$ denote the stationary subgroup of a point $x\in X.$

\begin{dfn} Suppose there is a subgroup $G_*(X)\subset G$ such that $G_x$ is a conjugate of $G_*(X)$ for $x$ in a dense subset in $X.$ Then we call $G_*(X)$ a
stabilizer in general position (SGP) for the action of $G$ on $X.$
\end{dfn}

For a reductive $G,$ and $X=V$ a vector space, a SGP $G_*(V)$ always exists \cite{Rich}, \cite{Lun}.

\begin{Lemma}(\cite{Oni})\label{Onishchik} Assume $H,G,G_*:=G_*(X)$  are reductive. Then
$(H,G,X)$ is exceptional if and only if the triple $(G,H,G_*)$ is a factorization, i.e., $G=HG_*.$
\end{Lemma}

\begin{prop}\label{centers} Assume $H,G,G_*$ are reductive. Then $(H,G,X)$ is exceptional if and only if $(G',H',(G_*)')$ and $(Z(G), Z(H)_{Z(G)}, Z(G_*)_{Z(G)})$ are factorizations, where $Z(G)$ denotes the center of $G,$ and $Z(H)_{Z(G)}, Z(G_*)_{Z(G)}$ denote the projections of the centers $Z(H)$ and $Z(G_*)$ into $Z(G).$
\end{prop}
\bpr This is a composition of Lemma~\ref{Onishchik} and Theorem 3.2 from \cite{Oni2}.\qed

\begin{Lemma}\label{vnizz} Assume $H,G,G_*$ are reductive. If $(H,G,X)$ is exceptional then
the triple $(H',G',X)$ is exceptional.
\end{Lemma}

However, the converse does not hold. Consider a triple $(H,G,X)=(\{id\}, \mathbb C^*,\mathbb C).$ Then $(H',G',X)=(\{id\},\{id\},\mathbb C)$ is exceptional, where $(H,G,X)$ is not.

\bpr
As follows from Proposition~\ref{centers},
$(G,H,G_*(X))$ is a factorization implies $(G',H',(G_*(X))')$ is a factorization. For  $x\in X,$ we
have $(G_x)'\subseteq G_x\cap G'=(G')_x.$ In particular,  $(G_*(X))'\subseteq (G')_*(X).$ Hence, $(G',H',(G')_*(X))$
is also a factorization, that is, $(H',G',X)$ is exceptional.\qed

\begin{Lemma}\label{orelse}
Let $(H,G,X)$ be an exceptional triple, $G=G_1\times G_2\supset H=H_1\times G_2.$ Then either
$(H_1,G_1,X)$ is exceptional, or $\dim pr_{G_2}G_*(X)>0.$
\end{Lemma}

\bpr Consider the maps
\[\xymatrix{
 {X/H_1} \ar[r]^f\ar[d]_{\pi_H}&{X/G_1}\ar[d]^{\pi_G}\\
{X/H}\ar@{=}[r]&{X/G,}\\
}\]
where $f$ is $G_2$-equivariant, and $\pi_H,$ $\pi_G$ are the morphisms of factorization by $G_2.$ Assume $\dim pr_{G_2}G_*(X)=0.$ Since $\dim pr_{G_2}G_*(X)=(G_2)_*(X/G_1),$ a generic $\pi_G$-fiber has (the maximal)
dimension equal to $\dim G_2.$ Since the diagram is commutative, and the dimension of a $\pi_H$-fiber is less or equal
$\dim G_2,$ we obtain $f=id,$ that is, $(H_1,G_1,X)$ is exceptional.\qed

\begin{claim}\label{normal}
Let $N\subset G$ be a normal subgroup. Then $\mathbb C(X)^G=\mathbb C(X/N)^{G/N}.$
\end{claim}

\section{Outlines}\label{out} In this section we proof Proposition~\ref{reduc} and outline the proof of Theorem~\ref{ortthm}.
\subsection{Proof of Proposition~\ref{reduc}.}

Let $H\subset G\subset SO(V)$ be orthogonal reductive linear groups.
Proposition~\ref{centers} allows us to assume that $Z(H)=Z(H)_{Z(G)},$ i.e., $Z(H)\subseteq Z(G).$ By definition of $T,$ we have $Z(G)=Z(H)\times T.$

 Assume $(H,G,V)$ is exceptional. Since a SGP of a reductive orthogonal group is reductive (\cite{Lun}, see also \cite{Vin0}), Lemma~\ref{vnizz} proves $(i).$ Furthermore, since $Z(H)$ has a trivial SGP on $V/(G'\times T),$
 the triple $(H',G'\times T,V)$ is exceptional by Lemma~\ref{orelse}, that is, $T$ acts trivially on $\mathbb C[V]^{G'}=\mathbb C[V]^{H'},$ and we get $(ii).$

Now suppose conditions $(i),(ii)$ hold.
Denote by $X:=V/G'=V/H'.$ By $(ii)$, we have $\mathbb C[X]^T=\mathbb C[X].$
By Claim~\ref{normal}, we have
$\mathbb C[V]^G=\mathbb C[X]^{Z(G)}=\mathbb C[X]^{T\times Z(H)}=\mathbb C[X]^{Z(H)}=\mathbb C[V]^H,$ i.e., $(H,G,V)$ is exceptional.
\qed

\subsection{Outlines of the proof of Theorem~\ref{ortthm}.}
We divided Theorem~\ref{ortthm} into Theorems~\ref{LTthm},~\ref{SFthm}, and~\ref{general}.
Theorems~\ref{LTthm} and~\ref{SFthm} treat two particular cases,
namely, locally trivial and  strongly faithful triples, and we prove them later in sections~\ref{LT} and~\ref{SF}. Theorem~\ref{general} is proved in this section.

As Lemma~\ref{Onishchik} suggests, one can try to apply
Onishchik's classification of factorizations for reductive groups \cite{Oni} to classify exceptional triples.  However, we cannot apply this
classification directly, since for an arbitrary complicated $G$-module $V$ one
cannot say much about $G_*(V),$ see an example at the end of section~\ref{diff}.  What we do is we apply Onishchik's classification
to  locally trivial and to strongly faithful triples,  and then
treat the general case manually based on these particular cases, see Theorem~\ref{general}.

\begin{prop}\label{diagonal} Suppose $(H,G,V)$ is semisimple orthogonal exceptional, $H$ is maximal in $G.$
$(H,G,V)$ is locally trivial if and only if $H$ is diagonal (see Definition~\ref{maximal}).
\end{prop}

\bpr The "only if" implication follows from the definition of a locally trivial triple.

Suppose $H$ is diagonal and assume $(H,G,V)$ is not locally trivial. Then there exists a non-zero
$G$-irreducible submodule $W\subset V$ such that $H_1\times H_2$ acts nontrivially on $W.$ Let
$W=\phi_1\otimes\phi_2\otimes\phi_{12},$ where $\phi_1,$ $\phi_2,$ $\phi_{12}$ are irreducible
representations of $H_1,$ $H_2,$ and $H_{12}$ respectively. By Lemma~\ref{spacedecom}, either $W$ is
$H$-irreducible, or $W=U\oplus U^*,$ where $U$ is $H$-irreducible. That is,  $\phi_1\otimes\phi_2$ restricted
to $H_{12}$, is equal either to (a) $\psi$ or to (b) $\psi\oplus\psi^*,$ where $\psi$ is
$H_{12}$-irreducible.

The option $(a)$ is impossible, see, for example, \cite{Kum}. Assume (b) holds. Denote by $\Lambda_1$ and
$\Lambda_2$ the highest weights of $\phi_1$ and $\phi_2.$ Then $\Lambda=\Lambda_1+\Lambda_2$ is the highest
weight of $\psi,$ and $\Lambda-M$ is the highest weight of $\psi^*,$ where $M$ is a non-negative linear
combination of simple roots of $H_{12}.$

In particular, we have $(\Lambda-M,\Lambda-M)=(\Lambda,\Lambda).$ However, this is impossible. We have
$(\Lambda-M,\Lambda-M)=(\Lambda,\Lambda)-(\Lambda,M)-(M,\Lambda-M)<(\Lambda,\Lambda),$ since $\Lambda-M$ is
dominant. \qed

\begin{thm}\label{LTthm} Suppose $(H,G,V)$ is semisimple indecomposable locally trivial orthogonal exceptional triple, $H$ is maximal in $G.$ Then $(H,G,V)$ has diagram one of  $\bf{T1-T9}.$
Conversely, all diagrams $\bf{T1-T9}$ correspond to exceptional triples, assuming that in diagram $\bf{T9}$ one of the two indecomposable components of $G$ satisfies Lemma~\ref{yak}.
\end{thm}

\begin{thm}\label{SFthm} Suppose $(H,G,V)$ is semisimple indecomposable strongly faithful orthogonal exceptional triple, $H$ is maximal in $G.$ Then $(H,G,V)$ has diagram one of $\bf{F1-F10}.$ Conversely, all diagrams
$\bf{F1-F10}$ correspond to exceptional triples.
\end{thm}

We prove these theorems in sections~\ref{LT} and~\ref{SF}.

Let us consider the general case.

\begin{thm}\label{general}
Assume a semisimple indecomposable orthogonal exceptional triple  $(H,G,V)$ is neither locally trivial, nor  strongly faithful. Then  $(H,G,V)$ has a diagram
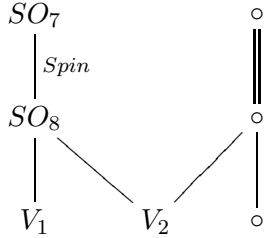
\[\xymatrix{
{SO_7} \ar@{-}[d]^{Spin}& &{\circ}\ar@{=}[d]\\
{SO_8} \ar@{-}[d]\ar@{-}[dr]& &{\circ} \ar@{-}[d]\ar@{-}[dl]\\
{V_1}&{V_2}&{\circ}\\
}\]
where the action of the subgroup $SO_8\subset G$ on $V_1\oplus V_2$ is isomorphic, up to a trivial component,
to $t\phi_1,$  $t\le 3.$ Conversely, any triple with such diagram   is exceptional.
\end{thm}

\bpr In this proof, all vertices in all diagrams correspond to simple subgroups/$G$-submodules, unless specified otherwise.

Since $(H,G,V)$ is not locally trivial, $H$ is a straight (see Definition~\ref{maximal}) subgroup in $G$ by
Proposition~\ref{diagonal}, that is, the top of a diagram for $(H,G)$ looks like
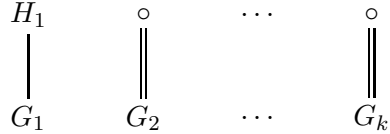
\[\xymatrix{
 {H_1} \ar@{-}[d]&{\circ}\ar@{=}[d]&{\dots} &{\circ}\ar@{=}[d]\\
{G_1}&{G_2}&{\dots}&{G_k}\\
}\]
where $H_1$ is maximal in $G_1.$ Since $G$ is not strongly faithful, we have $k>1.$ Assume $G_2\ne\{id\}.$
Denote by $\tilde V\subset V$ the sum of all submodules where $G_1$ acts nontrivially, that is, the sum of
all submodules connected to $G_1.$ Since $(H,G,V)$ is indecomposable, one of  $G_2,\dots, G_k,$ say $G_2,$ acts
nontrivially on $\tilde V,$ i.e., $G_2$ is connected to $\tilde V.$

Let $\tilde W\supseteq \tilde V$ be a simple submodule connected to $G_2,$ and $W\supset\tilde W$ be the
minimal orthogonal submodule containing $\tilde W.$ By Lemma~\ref{spacedecom}, either $W=\tilde W,$ or
$W=\tilde W\oplus\tilde W^*.$

Denote by $G^W$ the product of all normal subgroups in $G$ connected to $W$ (as well as to $\tilde W$), and
let $H^W$ be the projection of $G^W$ on $H.$  Then, by Corollary~\ref{sub}, $(H^W,G^W,W)$ is
exceptional, and, by definition of $W,$ strongly faithful. Since $G_1\subset G^W,$ we have $H^W\ne G^W,$ and
therefore, by Theorem~\ref{SFthm}, $(H^W,G^W,W)$ has a diagram

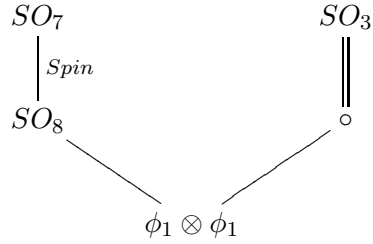
\[\xymatrix{
 {SO_7}\ar@{-}[d]^{Spin}& &{SO_3} \ar@{=}[d]\\
 {SO_8} \ar@{-}[dr]& &{\circ} \ar@{-}[dl]\\
 &{\phi_1\otimes \phi_1}& \\
}
\]

  This proves that
$(H,G,V)$ has a following diagram (here the vertices are not necessarily  simple):

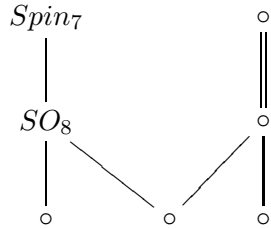
\[\xymatrix{
{Spin_7} \ar@{-}[d]& &{\circ}\ar@{=}[d]\\
{SO_8} \ar@{-}[d]\ar@{-}[dr]& &{\circ} \ar@{-}[d]\ar@{-}[dl]\\
{\circ}&{\circ}&{\circ}\\
}\]

Furthermore, Theorem~\ref{SFthm} implies $V=V^0\oplus V^1\oplus\dots\oplus V^k,$ where $(H,G,V^0)$ is a trivial triple, and  $(H,G,V^i),$ $i>0,$ is isomorphic to one of

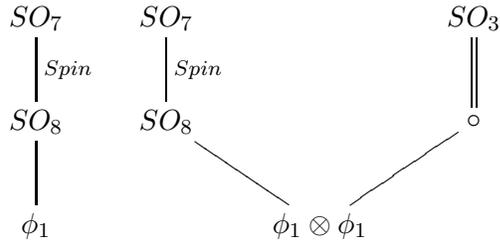
\[ \begin{tabular}{ll}
\xymatrix{
 {SO_7} \ar@{-}[d]^{Spin}\\
 {SO_8}\ar@{-}[d]\\
 {\phi_1}\\
}
&

\xymatrix{
 {SO_7}\ar@{-}[d]^{Spin}& &{SO_3} \ar@{=}[d]\\
 {SO_8} \ar@{-}[dr]& &{\circ} \ar@{-}[dl]\\
 &{\phi_1\otimes \phi_1}& \\
}
\\
\end{tabular}
\]
Hence,  $(H_1,G_1,V)$ is isomorphic, up to a
trivial component, to

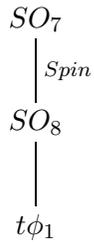
\[\xymatrix{
 {SO_7} \ar@{-}[d]^{Spin}\\
 {SO_8}\ar@{-}[d]\\
 {t\phi_1}\\
}
\]

 We only have to prove now that $t\le 3.$

\begin{Lemma}\label{t<=3}
$(H_1,G_1,V)$ is exceptional.
\end{Lemma}

By Theorem~\ref{SFthm}, Lemma~\ref{t<=3} implies $t\le 3.$

{\bf Proof of Lemma~\ref{t<=3}.} Denote $V^I=V^1\oplus\dots\oplus V^k,$ and let $G^I\subset G$ be the product
of all simple normal subgroups in $G$ acting nontrivially on $V^I.$ Let $H^I$ be the projection of $G^I$ on $H.$ By
Corollary~\ref{sub}, $(H^I,G^I,V^I)$ is exceptional.

We have $H^I=SO_7\times SO_3^{(1)}\times\dots\times SO_3^{(m)},$ $G^I=SO_8\times
SO_3^{(1)}\times\dots\times SO_3^{(m)}$ for some $m> 0$ (note that $m=0$ would mean $(H,G,V)$ is
decomposable).

Denote $S:=(SO_3^{(1)}\times\dots\times SO_3^{(m)})_*(V^I/SO_8).$ Let us show
that $\dim S~=~0.$

Note that each $SO_3^{(i)}$ is connected to exactly one of $V_1,\dots,V_k.$ Indeed, assume the opposite, say,
$SO_3^{(1)}$ is connected to $V_1$ and $V_2.$ Then $(H^I,G^I,V_1\oplus V_2)$ is isomorphic to

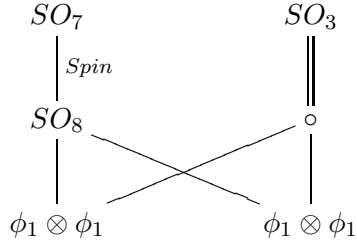
\[\xymatrix{
 {SO_7}\ar@{-}[d]^{Spin}& &{SO_3} \ar@{=}[d]\\
 {SO_8} \ar@{-}[d]\ar@{-}[drr]& &{\circ} \ar@{-}[d]\ar@{-}[dll]\\
 {\phi_1\otimes \phi_1}& &{\phi_1\otimes \phi_1} \\
}
\]

which is not exceptional by Theorem~\ref{SFthm}, in
contradiction to Corollary~\ref{sub}. Hence, we may assume for every $i=1\dots m$ that $SO_3^{(i)}$ is
connected to $V_i$ and disconnected from all other $V_j.$

Denote $S_i:=(SO_3^{(i)})_*(V_i/SO_8).$ Then $S_i=pr_{SO_3^{(i)}}(SO_8\times SO_3^{(i)})_*(V_i),$ and, therefore, $S_i$ contains $pr_{SO_3^{(i)}}S.$  However, $S_i\cong pr_{SO_3}(SO_8\times SO_3)_*(\phi_1\otimes~\phi_1)=\{id\},$ and therefore, $S=id.$  By
Lemma~\ref{orelse}, it means that $(H_1,G_1,V^I),$ as well as on $(H_1,G_1,V),$ is exceptional. Lemma~\ref{t<=3} is proved.\qed

Let us prove the last statement for the theorem.

Consider a triple $(H,G,V)$ that satisfies the conditions of the
theorem.  Then $(SO_7,SO_8,V)$ is exceptional by Theorem~\ref{SFthm}. In other words,
$V/SO_8=V/SO_7,$ and, by Claim~\ref{normal}, we have $\mathbb C[V]^G\cong \mathbb C[V/SO_8]^{G/SO_8}\cong\mathbb
C[V/SO_7]^{H/SO_7}\cong\mathbb C[V]^H.$ \qed

\section{Irreducible exceptional triples}\label{irrsummary}
This is a partial summary of \cite{irrex}.

Let $(H,G,V)$ be an irreducible semisimple triple.
Theorem~\ref{t1} classifies exceptional triples $(H,G,V)$ up to an explicit equivalence relation called castling
transform (Definition~\ref{defcst} below).

Suppose $F\subseteq SL(U)$.  Consider the group $G=F\times SL(W)\subset SL(U\otimes W),$ $\dim W\le \dim U$, and the group
$\check G=F\times SL(\check W)\subset SL(U^*\otimes \check W)$, where $\dim\check W=\dim U-\dim W.$
Then $\check G_*(U^*\otimes \check W)\cong G_*(U\otimes W)$ \cite{Ela2}.
In particular, if $\dim U=\dim W +1$, then $F_*(U)\cong G_*(U\otimes W).$

\begin{dfn}\label{defcst}We say that the (linear) group $\check G$ is an immediate castling transform of the group $G,$ and vice versa.
We say that a group $\tilde G$ is a castling transform of $G,$ and write $\tilde G\bowtie G,$ if $\tilde G$
is a result of a sequence of immediate castling transforms of $G$.
We will say that two triples $(\tilde H,\tilde G,\tilde V)$ and $(H,G, V)$ are congruent, and write $(\tilde H,\tilde
G,\tilde V)\bowtie (H,G,V)$ if they are isomorphic up to simultaneous castling transform.
\end{dfn}

One can show that if $(H,G,V)$ is  exceptional  and $(\tilde H,\tilde G,\tilde V)\bowtie (H,G,V),$ then $(\tilde H,\tilde
G,\tilde V)$ is also exceptional, and $\nu (\tilde H,\tilde G,\tilde V)=\nu (H,G,V).$

Denote $H_{s,t,k}=SL_s\times SL_t\times X_k\subset G_{s,t,k}=SL_{st}\times X_k\subset SL_{stk},$  where
$X_k\subseteq SL_k$ is irreducible, and $st>k.$ The group $H_{s,t,k}$ may have zero, one, or more invariants on $\mathbb C^{stk},$
depending on $s$,$t,$ and $k$ values, and also on $X_k,$ while $G_{s,t,k}$ has zero invariants.

\begin{thm}\label{t1}
Suppose $(H,G,V)$ is exceptional. If $\nu (H,G,V)=0,$ $(H,G,V)$ is congruent to $(H_{s,t,k}, G_{s,t,k},\mathbb C^{stk})$
or to one of the triples {\bf L1-L7}. Otherwise, $(H,G,V)$ is congruent to one of the triples {\bf A1-A10}.
Conversely, all triples {\bf L1-L7, A1-A10}, as well as their congruents, are exceptional.

\end{thm}

In diagram {\bf L5}, $Y_k$ is a maximal subgroup in $X_k\subseteq SL_k.$
For triples {\bf A1-A10}, we indicated the generators degrees for the
algebra of invariants, which is always polynomial.

\par\bigskip
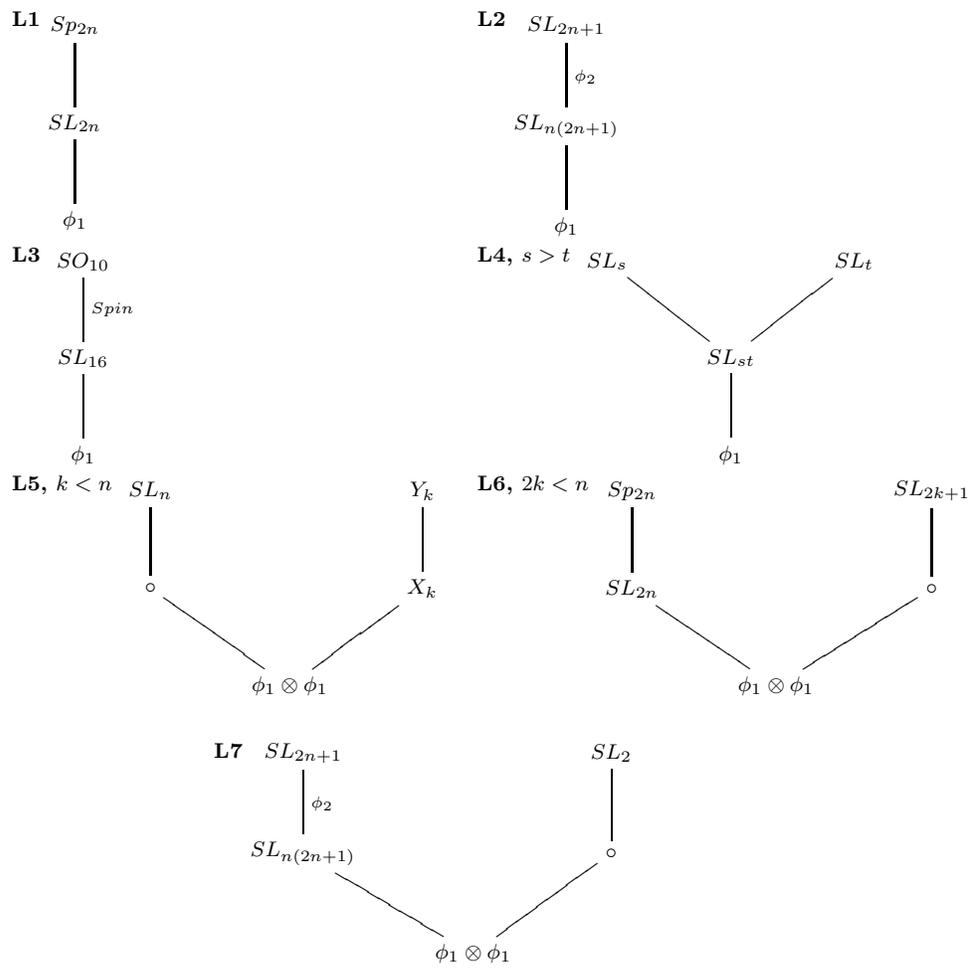
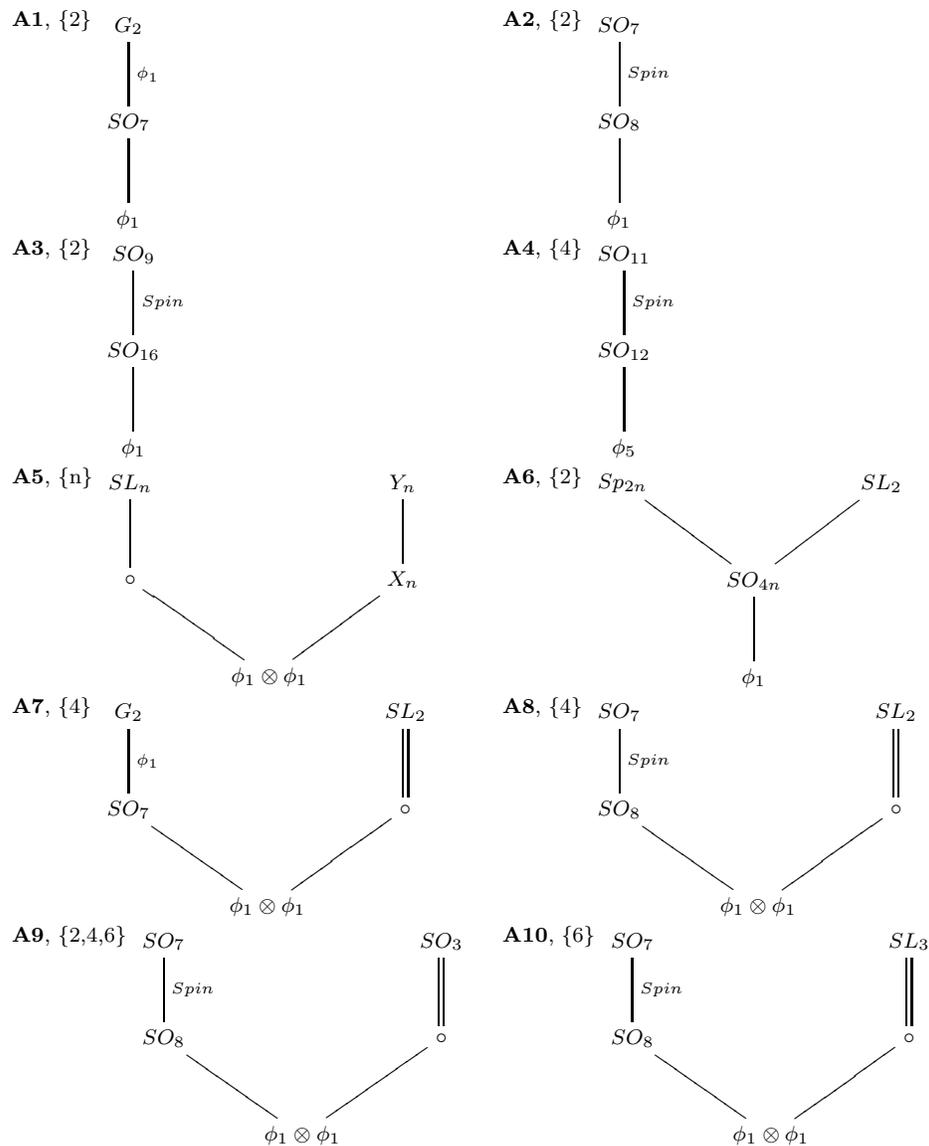
\begin{scriptsize}

\begin{tabular}{ll}
\bf{L1}\xymatrix{
 {Sp_{2n}} \ar@{-}[d]\\
 {SL_{2n}}\ar@{-}[d]\\
 {\phi_1}\\
}
&
\bf{L2}\xymatrix{
 {SL_{2n+1}} \ar@{-}[d]^{\phi_2}\\
 {SL_{n(2n+1)}}\ar@{-}[d]\\
 {\phi_1}\\
}
\\
\bf{L3} \xymatrix{
 {SO_{10}} \ar@{-}[d]^{Spin}\\
 {SL_{16}}\ar@{-}[d]\\
 {\phi_1}\\
}
&

\bf{L4}, $s>t$ \xymatrix{
 {SL_s} \ar@{-}[dr]& &{SL_t} \ar@{-}[dl]\\
 &{SL_{st}} \ar@{-}[d]& \\
 &{\phi_1}& \\
}
\\

\bf{L5}, $k<n$ \xymatrix{
 {SL_n} \ar@{-}[d]& &{Y_k} \ar@{-}[d]\\
 {\circ} \ar@{-}[dr]&& {X_k}\ar@{-}[dl]\\
 &{\phi_1\otimes\phi_1}& \\
}
&

\bf{L6}, $2k<n$ \xymatrix{
 {Sp_{2n}} \ar@{-}[d]& &{SL_{2k+1}} \ar@{-}[d]\\
 {SL_{2n}} \ar@{-}[dr]&& {\circ}\ar@{-}[dl]\\
 &{\phi_1\otimes\phi_1}& \\
}
\\
\end{tabular}
\[\bf{L7} \xymatrix{
 {SL_{2n+1}} \ar@{-}[d]^{\phi_2}& &{SL_2} \ar@{-}[d]\\
 {SL_{n(2n+1)}} \ar@{-}[dr]&& {\circ}\ar@{-}[dl]\\
 &{\phi_1\otimes\phi_1}& \\
}
\]

\begin{tabular}{ll}

{\bf A1}, \{2\} \xymatrix{
 {G_2} \ar@{-}[d]^{\phi_1}\\
 {SO_7}\ar@{-}[d]\\
 {\phi_1}\\
}
&

{\bf A2}, \{2\} \xymatrix{
 {SO_7} \ar@{-}[d]^{Spin}\\
 {SO_8}\ar@{-}[d]\\
 {\phi_1}\\
}
\\

{\bf A3}, \{2\} \xymatrix{
 {SO_9} \ar@{-}[d]^{Spin}\\
 {SO_{16}}\ar@{-}[d]\\
 {\phi_1}\\
}
&
{\bf A4}, \{4\} \xymatrix{
 {SO_{11}} \ar@{-}[d]^{Spin}\\
 {SO_{12}}\ar@{-}[d]\\
 {\phi_5}\\
}
\\

{\bf A5}, \{n\} \xymatrix{
 {SL_n} \ar@{-}[d]& &{Y_n} \ar@{-}[d]\\
 {\circ} \ar@{-}[dr]&& {X_n}\ar@{-}[dl]\\
 &{\phi_1\otimes\phi_1}& \\
}
&

{\bf A6},  \{2\} \xymatrix{
 {Sp_{2n}} \ar@{-}[dr]& &{SL_2} \ar@{-}[dl]\\
 &{SO_{4n}} \ar@{-}[d]& \\
 &{\phi_1}& \\
}
\\

{\bf A7}, \{4\} \xymatrix{
 {G_2}\ar@{-}[d]^{\phi_1}& &{SL_2} \ar@{=}[d]\\
 {SO_7} \ar@{-}[dr]& &{\circ} \ar@{-}[dl]\\
 &{\phi_1\otimes \phi_1}& \\
}
&
{\bf A8}, \{4\} \xymatrix{
 {SO_7}\ar@{-}[d]^{Spin}& &{SL_2} \ar@{=}[d]\\
 {SO_8} \ar@{-}[dr]& &{\circ} \ar@{-}[dl]\\
 &{\phi_1\otimes \phi_1}& \\
}
\\
{\bf A9}, \{2,4,6\} \xymatrix{
 {SO_7}\ar@{-}[d]^{Spin}& &{SO_3} \ar@{=}[d]\\
 {SO_8} \ar@{-}[dr]& &{\circ} \ar@{-}[dl]\\
 &{\phi_1\otimes \phi_1}& \\
}
&

{\bf A10}, \{6\} \xymatrix{
 {SO_7}\ar@{-}[d]^{Spin}& &{SL_3} \ar@{=}[d]\\
 {SO_8} \ar@{-}[dr]& &{\circ} \ar@{-}[dl]\\
 &{\phi_1\otimes \phi_1}& \\
}
\\

\end{tabular}

\end{scriptsize}

\section{Factorizations}\label{f}

In this section we bring down definitions and theorems  from Onishchik's
paper~\cite{Oni2} (Theorems 3.1-3.3, 4.3), mostly simplified for our needs, to be used in sections~\ref{LT} and~\ref{SF}. Note that here we use the word "triple" for a triple of groups.

Assume $G,$ $H,$ $S$ are reductive groups.

A triple of groups $(G, H, S),$ where $H,S\subset G,$ is called a {\it factorization} if $G=HS.$

$\g,$ $\h,\dots$ denote the tangent Lie algebras of the groups $G,$ $H,\dots.$ A triple of Lie algebras $(\g,
\h, \s),$ where $\h,\s\subset \g,$ is called a {\it factorization} if $\g=\h+\s.$

\begin{Lemma}\label{fac1}
A triple $(G,H,S)$ is a factorization if and only if a triple $(\g, \h, \s)$ is a factorization.
\end{Lemma}

Assume $\g$ and $\h$ are semisimple. For a reductive algebra $\a,$ denote by $\a'$ the semisimple part of
$\a.$

\begin{Lemma}\label{fac2}
A triple  $(\g, \h, \s)$ is a factorization if and only if a triple $(\g, \h, \s')$ is a factorization. In
particular, if $\s$ is commutative, then the factorization $(\g, \h, \s)$ is trivial.
\end{Lemma}

We say a triple $(\g, \h, \s)$ is {\it symmetric} to the triple $(\g,\s,\h).$ A triple symmetric to a factorization
is also a factorization.

\begin{Lemma}\label{fac1.5}
Suppose $\g$ is simple, $\h,\s$ semisimple,  and the triple  $(\g, \h, \s)$ is a nontrivial factorization. Then $(\g, \h, \s)$ is
isomorphic or symmetric to one of the factorizations from Table O. Conversely, any triple from Table O is a
factorization.
\end{Lemma}

\begin{tabular}{|l|l|l|l|l|} \multicolumn{4}{l}{\bf Table O.}  \\
\hline
$\g$&$\h$&$\h\hookrightarrow\g$&$\s$&$\s\hookrightarrow\g$\\
\hline
$\sl_{2n},$ $n>1$&$\sp_{2n}$&$\phi_1$&$\sl_{2n-1}$&$\phi_1\oplus id$\\
$\so_{2n},$ $n>1$&$\sl_{n}$&$\phi_1\oplus\phi_{n-1}$&$\so_{2n-1}$&$\phi_1\oplus id$\\
$\so_{4n},$ $n>1$&$\sp_{2n}\oplus\sl_2$&$\phi_1\otimes\phi_1$&$\so_{4n-1}$&$\phi_1\oplus id$\\
$\so_7$&$\g_2$&$\phi_1$&$\so_6$&$\phi_1\oplus id$\\
$\so_7$&$\g_2$&$\phi_1$&$\so_5$&$\phi_1\oplus 2id$\\
$\so_8$&$\so_7$&$\spin$&$\so_7$&$\phi_1\oplus id$\\
$\so_8$&$\so_7$&$\spin$&$\so_6$&$\phi_1\oplus 2id$\\
$\so_8$&$\so_7$&$\spin$&$\so_5$&$\phi_1\oplus 3id$\\
$\so_{16}$&$\so_9$&$\spin$&$\so_{15}$&$\phi_1\oplus id$\\
\hline
\end{tabular}
\par\medskip

A semisimple algebra is called {\it strongly semisimple} if all its ideals have rank $> 1.$ For a reductive $\a,$
we denote by $\a^s$ the maximal strongly semisimple subalgebra of $\a.$ We call $\a^s$ {\it the strongly semisimple
part of $\a.$}
We have a decomposition $\a=\a^s\oplus\a^r,$ where $\a^r$ is a sum of all rank 1 ideals of $\a$ and the
center $\z(\a).$

\begin{Lemma}\label{fac3}
Let $(\h^r)_{\g^r},$ $(\s^r)_{\g^r}$ denote the projections of $\h^r,$ $\s^r$ on $\g^r.$ Then $(\g, \h, \s)$
is a factorization if and only if both  $(\g^s, \h^s, \s^s)$ and $(\g^r, (\h^r)_{\g^r}, (\s^r)_{\g^r})$ are
factorizations.
\end{Lemma}

Suppose $\s,$ $\g$ are strongly semisimple, $\s=\s_1\oplus\dots\oplus\s_l,$ where all $\s_i$ are simple
ideals. Assume $\h\subsetneq\g$ is a maximal semisimple subalgebra.

\begin{Lemma}\label{fac4}
 Suppose $\h$ is a straight subalgebra, i.e., $\h=\h_1\oplus\g_2\subset\g_1\oplus\g_2,$ where $\g_1$ is
simple, and $\h_1$ is maximal in $\g_1.$ Denote by $\s_{i1}$ the projection of $\s_i$ on $\g_1.$ Then
$(\g,\h,\s)$ is  a factorization if and only if
 $(\g_1,\h_1,\s_{i1})$ is  a factorization for some $i.$
\end{Lemma}

\section{Locally trivial pairs}\label{LT}
In this section we prove Theorem~\ref{LTthm}.

Let $(H,G,V)$ be an orthogonal locally trivial triple, $H$ is maximal in $G.$ According to
Proposition~\ref{diagonal}, $H$ is a diagonal subgroup, that is, $H=H_{12}\times F\subset G=H_1\times
H_2\times F,$ where $H_{12}\cong H_1\cong H_2$ are simple, and $H_{12}$ is embedded diagonally in
$H_1\times H_2.$ Denote $G_*:=G_*(V).$ According to Lemma~\ref{Onishchik}, $(H,G,V)$ is exceptional if and only if $(G,H,G_*)$ is a factorization.
Note that, since $F\subset H,$ the triple of groups $(G,H,G_*)$ is a factorization if and only if $(H_1\times H_2, H_{12}, S)$ is a factorization, where $S$ denotes the projection $pr_{H_1\times H_2}G_*.$

Denote $V=V_1\oplus V_2\oplus V_{12},$ where $H_1$ acts trivially on $V_2$, $H_2$ acts trivially on $V_1$,
and $H_1\times H_2$ acts trivially on $V_{12}.$

\begin{Lemma}\label{vniz}$(H,G,V)$ is exceptional if and only if
$(H_{12}, H_1\times H_2, (V_1\oplus V_2)/F)$ is exceptional.
\end{Lemma}

\bpr By Claim~\ref{normal}, $\mathbb C[V]^H=\mathbb C[V]^G$ $\Longleftrightarrow$
$\mathbb C[V/F]^{H_{12}}=\mathbb C[V/F]^{H_1\times H_2}$  $\Longleftrightarrow$
$\mathbb C[(V_1\oplus V_2)/F]^{H_{12}}=\mathbb C[(V_1\oplus V_2)/F]^{H_1\times H_2}.$\qed

\begin{Lemma}\label{LT2}
Suppose we have a reductive triple $(H_{12}, H_1\times H_2, X)$  where  $X:=X_1\times X_2,$ so that the restrictions $(H_{12}, H_1\times H_2, X_i)=(H_{12}, H_i,X_i)$
are  trivial, $i=1,2.$ Suppose the SGP's $S_{ii}:=(H_{12})_*(X_i),$ $i=1,2,$ are reductive. Then $(H_{12}, H_1\times H_2,X)$ is exceptional if and only if $(H_{12}, S_{11},
S_{22} )$ is a factorization.
\end{Lemma}

\bpr We have $(H_1\times H_2)_*(X)\cong S_{11}\times S_{22},$ and  $S_{12}:=(H_{12})_*(X)\cong (S_{11})_*(X_2).$  The triple $(H_{12}, H_1\times
H_2,X)$ is exceptional $\Longleftrightarrow$ $\overline{O_X(H_{12})}=\overline{ O_X(H_1\times H_2)}$ $\Longleftrightarrow$ $\dim O_X(H_{12})=\dim O_X(H_1\times H_2)$
$\Longleftrightarrow$  $\dim H_{12}-\dim S_{12}=\dim (H_1\times H_2)-\dim
(S_{11}\times S_{22})$ $\Longleftrightarrow$ $\dim H_{12}-\dim S_{22}=\dim S_{11}-\dim S_{12}$
$\Longleftrightarrow$ $\dim O_{X_2}(H_{12})=\dim O_{X_2}(S_{11})$ $\Longleftrightarrow$ the triple
$(S_{11}, H_{12},X_2)$ is exceptional $\Longleftrightarrow$ $(H_{12}, S_{11}, S_{22} )$ is a
factorization. \qed

Back to the original triple $(H,G,V)$, denote $S_i=G_*(V_i),$   $S_{ii}=(H_{12})_*(V_i/F),$ $i=1,2.$ Then $S_{ii}$ is isomorphic to the projection of $S_i$ on $H_i.$ In particular,
$S_{ii}$ is reductive.

\begin{Lemma}\label{LT3} If $(H,G,V)$ is exceptional, then necessarily
$(H_{12}, S_{11}, S_{22})$ is a factorization. If $(V_1\oplus V_2)/F= V_1/F\times V_2/F,$ this
condition is also sufficient.
\end{Lemma}

\bpr $(H,G,V)$ is exceptional if and only if $(H_{12}, H_1~\times~H_2,(V_1\oplus V_2)/F)$ is exceptional (Lemma~\ref{vniz}). Since $(V_1~\oplus~V_2)/F\longrightarrow
V_1/F\times V_2/F$ is surjective and $G$-equivariant, $(H_{12}, H_1\times H_2,V_1/F\times V_2/F)$ is exceptional (Lemma~\ref{factor}). Now we can apply Lemma~\ref{LT2}.\qed

{\bf Proof of Theorem~\ref{LTthm}.} Suppose $(H,G,V)$ is exceptional.  Following Lemma~\ref{LT3},
we are looking for orthogonal actions of $H$ on $V_i,$ $i=1,2,$ such that $(H_{12}, S_{11}, S_{22})$ is a factorization.

Suppose  $(H_{12}, S_{11}, S_{22})$  is a nontrivial factorization. Then we are in the conditions of
Lemma~\ref{fac1.5}. Therefore,  $(H_{12}, S_{11}, S_{22})$ is isomorphic or symmetric to one of the triples from Table O. Hence, we have to find all the actions $H:V_i$ such that $H_{12}$ and $S_{ii}$
appear in the same line in Table O. Let's call such actions candidates.

Note that $V$ may only have solid or split summands (see Definition~\ref{I-III}).  Using \cite{Ela1}, \cite{Ela2},
we first find all solid candidates. The result is  all the entries
in Table T except for the last one.

Now let us  prove that the only split candidate is the
last entry of Table T.
Let $V_i=U\oplus U^*$, where  the action $H:U$ is not orthogonal.

Assume $F$ is trivial. Then \cite{Ela1}, \cite{Ela2} imply that $H=H_{12}:U$ is isomorphic to $SL_{2n}:\phi_1\oplus\phi_1^*.$

Let us prove that there does not exist a split candidate with $F$ not trivial. We may assume that $U$ is irreducible.

Note that an action $H:V_i$ is a candidate $\Longleftrightarrow$ $(S_{ii}, H_{12}, V_i)$ is exceptional $\Longrightarrow$ $(S_{ii}, H_{12}, U)$ is exceptional.

Consider  $H$ acting on $U$ as $SL_n\times SL_k:\phi_1\otimes\phi_1,$ $n\ge k\ge 2.$ Then $S_{ii}=SL_{n-k}\times T_{k-1},$ where $T_{k-1}$ is a torus. Since
$k\ge 2,$
$(S_{ii}, H_{12},U)$ is not exceptional. Even more so $(S_{ii}, H_{12},U)$ is not exceptional when  $H$ acts on $U$ as $SL_n\times X_k:\phi_1\otimes\phi,$ $X_k\subseteq SL_k,$ whether $H_{12}=X_k$ or $H_{12}= SL_n.$
 Note that since only the groups $SL_n$ and
$SO_n$ appear in the first column of Table O, the last option we have for a candidate is   $H_{12}=SO_n,$ $F\subseteq SL_k,$ $U=\phi_1\otimes\phi.$  Here we have $S_{ii}=SO_{n-2k},$ or $id$ if $n\le 2k.$ Since $k\ge 2,$ this action is not a candidate.

We proved that $H:V_i$ has to be isomorphic to one of the following:
\par\medskip
\begin{tabular}{|l|l|l|l|} \multicolumn{4}{l}{\bf Table T.}  \\
\hline
$H$& $H:V_i$&$S_{ii}$&$S_{ii}\hookrightarrow H$\\
\hline
$SL_4$&$\phi_2$&$Sp_4$&$\phi_1$\\
$SO_7$&$Spin$&$G_2$&$\phi_1$ \\
$SO_{12}$&$Spin$&$SL_6$&$\phi_1\oplus\phi_5$ \\
$SO_7$&$\phi_1$&$SO_6$&$\phi_1\oplus id$\\
$SO_8$&$Spin$&$SO_7$&$\phi_1$ \\
$SO_{2n}$&$\phi_1$&$SO_{2n-1}$&$\phi_1\oplus id$\\
$SO_7$&$2\phi_1$&$SO_5$&$\phi_1\oplus 2id$\\
$SO_8$&$2\phi_1$&$SO_6$&$\phi_1\oplus 2id$\\
$SO_8$&$3\phi_1$&$SO_5$&$\phi_1\oplus 3id$\\
$SO_8\times SO_3$&$\phi_1\otimes\phi_1$&$SO_5$& $(\phi_1\oplus 3id)\otimes id$\\
\hline
$SL_{2n}$&$\phi_1\oplus\phi_{2n-1}$&$SL_{2n-1}$&$\phi_1\oplus id$\\
\hline
\end{tabular}
\par\bigskip

It's now easy to make sums $V=V_1\oplus V_2\oplus V_{12}$ so that $(H,G,V)$ is exceptional.
The result is diagrams {\bf T1-T8}.

Now suppose  $(H_{12}, S_{11}, S_{22})$  is a trivial factorization, say, $H_1=S_{11}.$
As in the beginning of the section, let $S$ denote the projection of $G_*$ to $H_1\times H_2.$

\begin{Lemma}\label{sl2}
Suppose $H_1=S_{11}.$ Then

(a) $H_1=SL_2.$ Moreover, for any minimal orthogonal subspace  $W\subseteq V_1,$ the action
 $H:W$ is isomorphic to  $SL_2\times Sp_{2n}:\phi_1\otimes\phi_1,$ $n\ge 1.$

(b) $(H,G,V)$ is exceptional if and only if $S\ne H_{12}.$
\end{Lemma}

\bpr
(a) follows from \cite{Ela1}, \cite{Ela2}.

We have $(G,H,G_*(V))$ is a factorization $\Longleftrightarrow$ $(SL_2\times SL_2, SL_2^{diag}, S)$ is a factorization. Since the projection of $S$ on the first factor is $SL_2,$   this triple is always a factorization unless $S=SL_2^{diag}=H_{12}.$ Lemma~\ref{sl2} is proved.\qed

In \cite{yak}, we find the following
description of  groups $K=SL_2\times\check K\subset SO(V)$ such that the projection of SGP $K_*:=K_*(V)$ on the first factor is $SL_2.$

Consider a rooted tree $T$ with vertices $0,1,\dots, q,$ where $0$ is the root. Assign a weight $d(i)$ to each vertex such that  $d(i)$ is either a positive integer or $\infty.$ Assume that (a) $d(0)=1,$ (b) each vertex $i$ with $d(i)=\infty$ has degree $1$ (here the degree of a vertex is the number of adjacent vertices), and (c) if $(i,j)$ is an edge with $d(j)=\infty,$ then $d(i)>1.$ We say a vertex $i$ is {\it finite} if $d(i)<\infty,$ and an edge $(i,j)$ is {\it finite} if both $i$ and $j$ are finite.

Let $\check K$ be a product of $Sp_{2d(i)}$ over all finite vertices except for the root, and $K$ be a product of $Sp_{2d(i)}$ over all finite vertices. Assign a vector space $W_{(i,j)}$  to each edge $(i,j),$ namely,  let  $W_{(i,j)}=\mathbb C^{2d(i)}\otimes\mathbb C^{2d(j)}$ for a finite $(i,j),$  and  $W_{(i,j)}=S_0^2\mathbb C^{2d(i)}$  for an edge $(i,j)$ with $d(i)=\infty.$

Now let $V$ be a direct sum of $W_{(i,j)}$ over all edges of $T.$ The group $K$ acts on $V$ in the following way: each factor $Sp_{2d(i)}$ of $K$ acts (a) on the first factor of $W_{(i,j)}$ for all finite vertices $j$ connected to $i$ by an edge, and also (b) on $W_{(i,j)}$ for all infinite $j$ connected to $i$ by an edge. For example, a tree with $2$ vertices corresponds to the group $K=SL_2\times Sp_{2d(1)}$ acting on $\mathbb C^2\otimes\mathbb C^{2d(1)}.$

\begin{Lemma}\label{yak}(\cite{yak}) Consider a tree as above and assume that the corresponding linear group $K=SL_2\times\check K\subset SO(V)$ is indecomposable, i.e., $K$ is not a direct sum of proper linear subgroups. The projection of $K_*:=K_*(V)$ on the first factor is $SL_2$ if and only if

(I) all vertices $i$ with $d(i)>1$ have degree at most 2;

(II) for all edges $(i,j)$ such that $d(i)>1,$ $d(j)>1,$ one of the vertices $i,j$ has degree 1.
\end{Lemma}

Suppose a tree $T$ satisfies the conditions of the lemma  and assume that $d(1)=1.$  It's easy to check that if vertex $1$ belongs to the connected component of the root, then the projection of $K_*$ on the subgroup $Sp_{2d(0)}\times Sp_{2d(1)}\cong SL_2\times SL_2$ is the diagonal subgroup $SL_2^{diag}.$ If $0$ and $1$ are disconnected, then this projection is  $SL_2\times SL_2.$

Consider the tree $T_G$ corresponding to $G.$ Since $(H,G,V)$ is indecomposable, $T_G$ is either connected or $T_G$ has two connected components: one for $H_1$, and one for $H_2.$ We saw that $S=SL_2^{diag}$ if and only if  $H_1$ and $H_2$ belong to one connected component on $T_G.$ Hence, by Lemma~\ref{sl2}, we obtain the following diagram {\bf T9} corresponding to an exceptional $(H,G,V):$

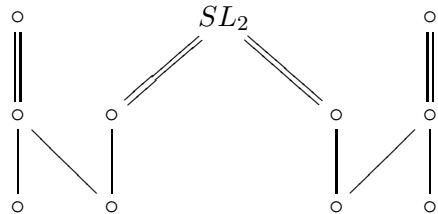
\[
\xymatrix{
 {\circ} \ar@{=}[d]&  &{SL_2}\ar@{=}[dl]\ar@{=}[dr]& &{\circ} \ar@{=}[d]\\
{\circ} \ar@{-}[dr]\ar@{-}[d] &{\circ} \ar@{-}[d]& &{\circ}\ar@{-}[d]&{\circ}\ar@{-}[dl]\ar@{-}[d]\\
{\circ}&{\circ}&  & {\circ}& {\circ}\\
}
\]

In this diagram, we replaced all finite vertices of $T_G$ from the connected components of $H_1$ and $H_2$ by the very left and very right circles in the middle level of the diagram correspondingly, and  then substituted all the edges with  vertices in the bottom level.

The theorem is proved.\qed

\section{Strongly faithful pairs}\label{SF}
Here we prove Theorem~\ref{SFthm}.

Let $(H,G,V)$ be an orthogonal strongly faithful triple, $H$ is maximal in $G.$ By
Proposition~\ref{diagonal}, $H$ is a straight subgroup: $H=H_1\times G_2\subset G=G_1\times
G_2,$ where  $H_1$ is maximal in
 $G_1,$ and $G_1$ is simple.

 Take $V_1\subset V$ a minimal orthogonal $G$-invariant
subspace.

Suppose $V_1$ is solid (see Definition~\ref{I-III}). Then, by Theorem~\ref{t1}, $(H,G,V_1)$ has
 diagram one of {\bf F1-F5}.

Suppose $V_1$ is half-split. Then, as follows from \cite{Ela1}, \cite{Ela2}, and Lemma~\ref{Onishchik},  $(H,G,V_1)$ has diagram {\bf F6}.

Now suppose $V_1$ is split, that is, $V_1=U\oplus U^*,$
where $(H,G,U)$ is irreducible exceptional
non-orthogonal. By Theorem~\ref{t1}, the action of $G$ on $U$ is
isomorphic to one of the actions in Table F1 below:
\par\medskip
\begin{tabular}{|l|l|l|} \multicolumn{3}{l}{\bf Table F1.}  \\
\hline
&$G$&$U$\\
\hline
1&$SL_n$&$\phi_1$\\
2&$SL_n\times X_k,$ $X_k\subseteq SL_k,$ $k\le n$&$\phi_1\otimes\phi$\\
3&$SL_{n_1}\times\dots\times SL_{n_s}\times X_k,$ $k< n_1\dots n_s$&$\phi_1\otimes\dots\phi_1\otimes\phi$\\
4&$SO_{12}$&$\phi_5$\\
5&$SO_7\times SL_k,$ $k=2,5$&$\phi_1\otimes\phi_1$\\
6&$SO_8\times SL_k,$ $k=2,3,5,6$&$\phi_1\otimes\phi_1$\\
\hline
\end{tabular}
\par\medskip

Following Lemma~\ref{fac4}, we are interested in projections
of  the SGP $G_*(V_1)$ onto simple normal
subgroups of $G.$  Note that in case 3 we are only interested in
the projections onto the last factor $X_k.$

In case 1, $G=SL_n,$ we obtain an exceptional (by Lemma~\ref{fac4}) triple with diagram {\bf F7}.

In case 2, take $X_k = SL_k.$ Then $G_*(V_1)= SL_{n-k}\times T_{k-1},$
where $T_{k-1}$ is a torus, and  $SL_{n-k}\subset SL_n.$ Since
$k>1,$ $(G,H,G_*(V_1))$ is not a factorization by Lemma~\ref{fac1.5},
that is, $(H,G,V_1)$ is not exceptional. This implies $(H,G,V_1)$ is not
exceptional for any $X_k\subset SL_k,$ and therefore, it is not
exceptional in case 3 either.

In cases 4, 5 and 6, $G_*(V_1)=SL_4,$ $SO_{7-2k}\subset SO_7,$ and
$SO_{8-2k}\subset SO_8,$ respectively,  and $(H,G,V_1)$ is not
exceptional by Lemma~\ref{fac1.5}.

Suppose there are other minimal orthogonal $G$-invariant subspaces
$V_2,$ $\dots,$ $V_k~\subset~V.$ From the explicit description we
obtained for the action of $G$ on $V_1,$ we see that necessarily
for all $i$ $V_i\cong V_1$ as a $G$-module. Then we simply check
directly  if $(H,G,V_1\oplus V_2\oplus\dots)$ is exceptional. Table F2 below shows
this calculation. As a result, we obtain diagrams {\bf F8-F10}  (cases 1,3,5 resp. in Table F2).

\par\medskip
\begin{tabular}{|l|l|l|l|} \multicolumn{4}{l}{\bf Table F2.} \\
\hline
&$G$&$G:V$&$G_*(V)$\\
\hline
1& $SO_7$&$2\phi_1$&$SO_5$ \\
2& $SO_7$&$3\phi_1$&$SO_4$ \\
3& $SO_8$&$2\phi_1$&$SO_6$ \\
4& $SO_8$&$\phi_1\oplus\phi_3$&$G_2$ \\
5& $SO_8$&$3\phi_1$&$SO_5$ \\
6& $SO_8$&$4\phi_1$&$SO_4$ \\
7& $SO_{2n}$&$2\phi_1$&$SO_{2n-2}$\\
8& $SO_8\times SO_3$&$2(\phi_1\otimes\phi_1)$&$T_1$\\
9& $SO_8\otimes SO_3$&$\phi_1\otimes\phi_1\oplus(\phi_1\otimes\phi_1)^*$&$id$\\
10&$SL_n$&$2(\phi_1\oplus \phi_{n-1})$&$SL_{n-2}$\\
\hline
\end{tabular}
\par\medskip
The theorem is proved.

\end{document}